\documentclass[12pt,a4paper]{amsart}

\usepackage[utf8]{inputenc}
\usepackage[T1]{fontenc}
\usepackage{amsmath,amsfonts,amssymb}
\usepackage{newtxtext,newtxmath}

\usepackage{hyperref}
\usepackage{graphicx}
\usepackage{enumerate}

\usepackage{caption,subcaption}
\usepackage{array}

\usepackage[all]{xypic}
\usepackage[margin=25mm,top=30mm]{geometry}

\renewcommand{\leq}{\leqslant}
\renewcommand{\geq}{\geqslant}

\usepackage{xcolor}
\usepackage{algorithmicx}
\usepackage{algpseudocode}
\usepackage{listings} 

\usepackage{amsrefs}
\usepackage{float}

\newtheorem{definition}{Definition}
\newtheorem{proposition}{Proposition}

\title{A New Non-Archimedean Metric on Persistent Homology}

\author{İSMAİL GÜZEL}
\email{iguzel@itu.edu.tr}

\author{ATABEY KAYGUN}
\email{kaygun@itu.edu.tr}
\address{Department of Mathematics, Istanbul Technical University, Istanbul, Turkey.}

\begin{document}

\maketitle

\begin{abstract}
  In this article, we define a new non-archimedean metric structure, called \emph{cophenetic
    metric}, on persistent homology classes of all degrees. We then show that zeroth
  persistent homology together with the cophenetic metric and hierarchical clustering
  algorithms with a number of different metrics do deliver statistically verifiable
  commensurate topological information based on experimental results we obtained on
  different datasets. We also observe that the resulting clusters coming from cophenetic
  distance do shine in terms of different evaluation measures such as silhouette score and
  the Rand index. Moreover, since the cophenetic metric is defined for all homology degrees,
  one can now display the inter-relations of persistent homology classes in all degrees via
  rooted trees.
\end{abstract}

\section{Introduction}
\label{intro}
In this article, we define a new non-archimedean metric (a.k.a. an ultra-metric) called
\emph{cophenetic metric} on persistent homology classes of all degrees using only
homological information. Then we statistically verify that the topological information
coming from the zeroth persistent homology with our cophenetic metric is consistent with the
information provided by different hierarchical clustering algorithms using different metrics
based on numerical experiments on different datasets. We also observe that the clusters we
obtained via the cophenetic metric do yield competitive silhouette scores and the Rand
indices in comparison with clusters obtained from other metrics.

  Persistent homology is a new class of topological invariants that has
  found applications in statistical data analysis. Its computations rely on techniques from
  topology, algebraic topology, computational geometry, and computational linear algebra.
  In simplest terms, this class of invariants (recorded as a basis of a vector space) aim to
  keep a record of the topological features of the manifold that our data is sampled from
  such as the connected components, circle-like and higher dimensional sphere-like
  cavities. Moreover, the calculated topological invariants of the manifold from which our
  data is sampled come with a scale or a time parameter. This allows homology classes we
  computed to evolve as the parameter changes (creation, annihilation, merging, splitting
  etc.) and yields opportunities for heuristic fine-tuning depending on needs and
  computational constraints.  Furthermore, while the scale parameter may interact with the
  ambient metric of the data space, since homological invariants are purely topological the
  persistent homology is impervious to any perturbations in the metric structure as well as
  the number of features of the data space. This may be useful in situations where the
  ambient metric needs to be adjusted depending on the constraints, or worse yet, where no
  such metric exists.

 In the standard representation of persistent homology, one only keeps a record of
the number of persistent homological classes in the form of
\emph{life-time intervals}, called \emph{barcodes}~\cite{carlsson2005persistencebarcode}.  However, persistent
homology classes carry a very rich combinatorial structure, and one can do more that just
counting them. In reference to the dendrograms of hierarchical clustering
  algorithms, Carlsson expresses the same idea as a question
in~\cite[Ch.8]{miller2020handbook} and \cite{carlsson2020persistent}:
\begin{quote}
	The dendrogram can be regarded as the “right” version of the invariant $\pi_0$ in the
	statistical world of finite metric spaces. The question now becomes if there are similar
	invariants that can capture the notions of higher homotopy groups or homology groups.
\end{quote}
Since hierarchical clustering schemes, and therefore, dendrograms and non-archimedean
metrics are known to be equivalent by Carlsson~\cite{carlsson2010hierarchical}, to answer this
question it is clear that one needs to define a non-archimedean metric for homology or
homotopy classes.

\paragraph{Our contributions.}
The main contribution of this paper is a new non-archimedean metric defined on persistent
homology classes in all degrees using purely homological information coming from the
changing scale parameter.  For this metric, we analyze how persistent homology classes of a
certain degree ``merge'' on top of recording the birth and death times of these classes as
the scale parameter changes. We observe that since all data points
  naturally appear as zeroth degree persistent homology classes, one can compare our metric
  with standard distance measures with respect to their performances in machine learning
  algorithms that rely on distance measures on data points.  We tested the soundness of our
  proposal by checking whether hierarchical clustering methods with different metrics do
  indeed yield statistically verifiable commensurate topological information. We observe
  that when we measure the clusters obtained from the cophenetic metric against the clusters
  coming from other metrics we obtained competitive results. Furthermore, the metric we
  define on the persistent homology can now be used to sketch rooted tree presentations of
  the persistent homology classes in all degrees that track how these classes \emph{merge} as
  the scale parameter changes.

In forming the bridge between hierarchical clustering and persistent homology, we also
found that the answer to the question raised by Carlsson~\cite{carlsson2020persistent}
comes from algebraic topology: cobordisms. Dendrograms are 1-dimensional cobordism classes
of disjoint union of points.  For higher homology classes, one has to resort to
$n+1$-dimensional cobordisms of disjoint unions of $n$-spheres.  For example, for
persistent homology in degree 1 such cobordisms are given by oriented genus-$g$ Riemann
surfaces with finitely many punctures, and the classification of such 2-manifolds is
complete. Unfortunately, in dimensions 2- and higher such cobordisms are very difficult to
classify.  We are going to investigate the special case of persistent homology of degree-1
in an upcoming paper.

\paragraph{Prior art.}
Topological data analysis (TDA) is a new data analysis discipline whose fundamentals
straddle both very abstract and concrete sub-disciplines of the mathematical research. Even
though the theoretical roots TDA are firmly placed in algebraic topology, to solve its
computational needs it heavily uses computational geometry and numerical linear algebra.
Since TDA relies on the topology rather than a particular metric
  structure of the ambient space from which data is sampled, in theory, it is more suitable
  for extracting information from data for which a canonical metric is not clear from the
  context, or worse yet, does not exist.

Clustering algorithms, on the other hand, have been around for a long time and they
form an important and well-understood class of machine learning
algorithms~\cite{lance1967general,cophenetic1973numerical,JainDubes88,Legendre2012}. They are known to be equivalent to non-archimedean
metrics~\cite{carlsson2010hierarchical}. For a given data set, these algorithms aim to
deliver an optimal partition where subsets are supposed to show a high degree of
heterogeneity between, and a high degree of homogeneity within each subset.  However,
one has to make unavoidable ad-hoc choices to determine an optimal hierarchical
clustering algorithm for any data set at hand due to the impossibility result of Kleinberg~\cite{kleinberg:impossibility}.

Similar to clustering algorithms, the TDA methods we investigate in this paper also
rely on a changing scale parameter. But instead of relying on the metric structure
alone, these methods propose using \emph{persistent homology} to compute topological
invariants of a data set. Persistent homology was first introduced to investigate
topological simplifications of alpha shapes in~\cite{edelsbrunner2000topological}, but
later extended to arbitrary dimensional spaces in~\cite{zomorodian2005computing}.  The
topological invariants that persistent homology identifies are \emph{the Betti numbers}
defined for every natural number $n$. For instance, the Betti numbers for $n=0, 1$ and
$2$ indicate respectively the number of connected components, 2- and 3- dimensional
holes within the data set.  The information that persistent homology yields on the
change in topological features as the filtration scale parameter increases can be
presented in various different ways by
barcodes~\cite{carlsson2005persistencebarcode} and \cite{ghrist2008barcodes}, by persistence
diagrams~\cite{cohen2007stability}, by landscapes~\cite{bubenik2015statistical}, by
images~\cite{image2017persistence}, by terraces~\cite{terrace2018persistence}, by
entropy~\cite{merelli2015topological} and by curves~\cite{chung2019persistence}.  Barcodes
are the most commonly used representations of persistent homology classes in which one
keeps a record of finite collections of scale parameter intervals over which individual
persistent homology classes \emph{persist}.

Hierarchical clustering algorithms that we consider in this paper extract their results
based solely on the metric structure of the ambient space where the data set is
embedded. One can also statistically test the stability and convergence of these
methods~\cite{carlsson2008persistent,carlsson2010hierarchical}. In addition, they use
a convenient tree representations, called \emph{dendrograms}, to display the
information on how these clusters merge as the underlying scale parameter
changes~\cite{johnson1967hierarchical}. One can compare dendrograms coming from
different variants of clustering algorithms using a suitable
metric~\cite{JardineSibson71,Hartigan85,JainDubes88}, or even improve forecasts by
comparing dendrogram-like features obtained from different hierarchical clustering
methods on the same data set~\cite{ignacio2020intrinsic}.

In \cite{elkin2020mergegram}, the authors develop diagrams called \emph{mergegrams} derived out of the dendrograms of the hierarchical clustering algorithms as a
  replacement for the zeroth bar code.  However, unlike the mergegrams that depend on the
  metric structure of the ambient space and work with the zeroth bar codes only, our
non-archimedean metric relies on purely homological information and are defined for all
persistent homology degrees.  More importantly, they can be used for purposes other than
sketching rooted tree presentations.

\paragraph{Plan of the article.}
This paper is organized as follows. We give the necessary background material we need
on persistent homology and hierarchical clustering in Sections \ref{sect:PH} and
\ref{sect:HC}. We outline our methodology in Section~\ref{sect:PHandHC}.  The results
of our numerical experiments\footnote{The source code and the data of the numerical
	experiments we conducted in the paper can be found on the authors' GitHub page at
	https://github.com/ismailguzel/TDA-HC .}  are given in Section
\ref{sect:experiments}, and in Section \ref{sect:conc_future} we present our detailed
analysis of our experimental work in the light of theoretical discussions we presented
in the earlier sections. We also propose several avenues of future work in the same
Section.

\section{Persistent Homology}\label{sect:PH}

\subsection{Point clouds and simplicial complexes}

In hierarchical clustering, a collection of points given in an ambient metric space carries
no information other than the distances between them.  Derived information such as
cophenetic matrices also rely on this metric structure.  However, there are other tools to
derive more information about the topology of the data set at hand.  One of these useful
tools one can use is a \emph{simplicial complex}.

An abstract simplicial complex $K$ in a space $X$ is a collection of subsets of $X$ such
that for any two $x,y\in K$ one also has $x\cap y\in K$.  There are two variants of
simplicial complexes that are of interest for us: Vietoris-Rips complexes and \v{C}ech
complexes.

\subsubsection{Vietoris-Rips complexes}

Given a point cloud $D$, the Vietoris-Rips is defined to be the simplicial complex whose
simplicies are all points in $D$ that are at most $ \varepsilon $ apart.
\[ R_{\varepsilon}(D) = \{\sigma \subset D\mid d(x,y)\leq \varepsilon, \text{ for all }
x,y\in \sigma\} \]

\subsubsection{\v{C}ech complexes} 

Given a point cloud $ D $, the \v{C}ech complex associated with $D$ is defined to be the
simplicial complex given by
\[\mathcal{C}_{\varepsilon}(D) = \left\{\sigma \subseteq D \mid
\bigcap_{x \in \sigma} B_{\varepsilon}(x) \neq \varnothing\right\}. \] In other words, a
collection of points $\sigma = (x_0,\ldots,x_\ell)$ forms an $\ell$-simplex if the set of
balls of radius $ \varepsilon $ centered at these points has non-empty intersection.

\subsection{Choosing an appropriate scale parameter}

\begin{figure}[ht]
	\begin{subfigure}{.3\textwidth}
		\centering
		\includegraphics[scale=.2]{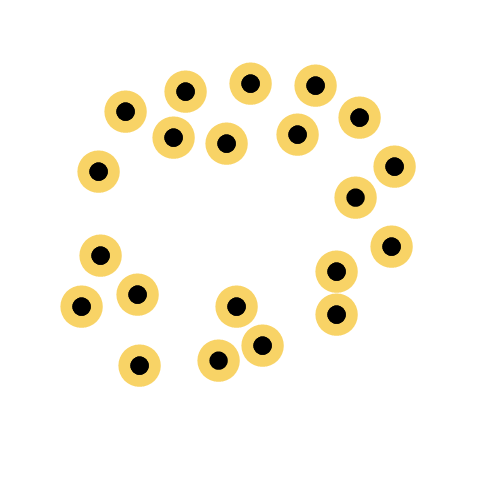}
		\caption*{$R_{\varepsilon_1}$}
	\end{subfigure}
	\begin{subfigure}{.3\textwidth}
		\centering
		\includegraphics[scale=.2]{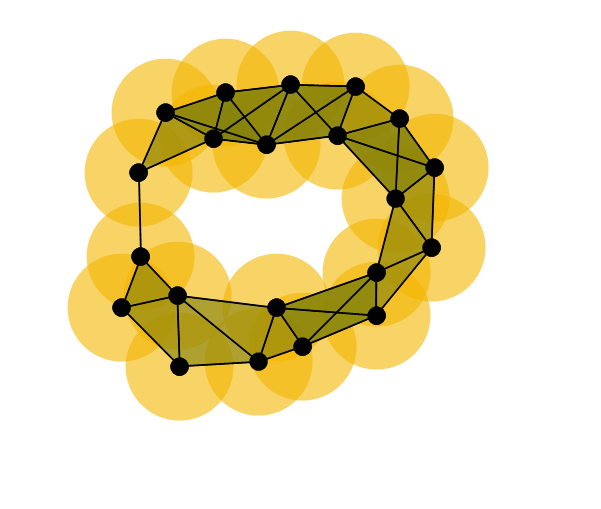}
		\caption*{$R_{\varepsilon_2}$}
	\end{subfigure}
	\begin{subfigure}{.3\textwidth}
		\centering
		\includegraphics[scale=.2]{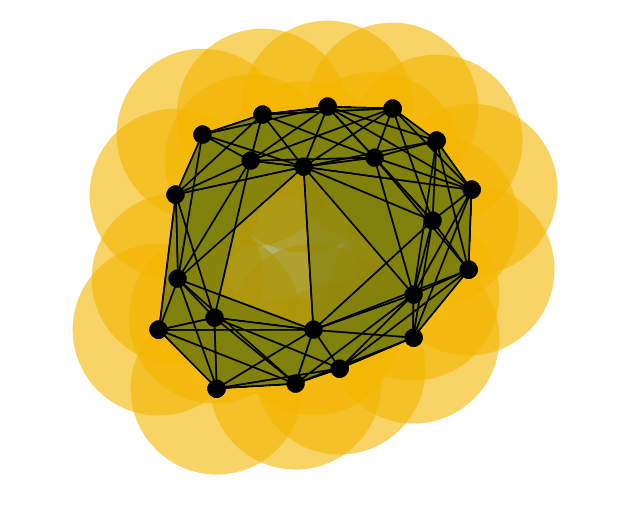}
		\caption*{$R_{\varepsilon_3}$}
	\end{subfigure}
	\caption{Vietoris-Rips complexes with increasing values of the parameters.}\label{epsilonincreasing}
\end{figure}

In order to turn a point cloud $D$ into a simplicial complex, we are going to use
$R_\varepsilon(D)$ the Vietoris-Rips complex associated with $D$ with a chosen proximity
parameter $ \varepsilon>0 $.  We then try to capture the topological features of the data by
changing the parameter $\varepsilon$.  As we see in Figure \ref{epsilonincreasing}, we may
not capture all of the topological features of data for a given proximity $ \varepsilon$.
Finding the optimal value for $\varepsilon$ for a given data set $D$ is a challenging
problem.

In \cite{edelsbrunner2000topological} and~\cite{zomorodian2005computing}, the authors proposed that the persistence homology might help to
determine an optimal value for $ \varepsilon $.  In persistent homology, one records the
longevity of each topological feature (in this case homology classes of $R_\varepsilon(D)$)
of a given data set as the proximity parameter $ \varepsilon $ changes.  One does this by
observing the \emph{persistence} of these topological features depending on $\varepsilon$.

\subsection{Persistent Homology}

Let $ \{ K_{\varepsilon}\mid \varepsilon \in \mathbb{R}_{+} \} $ be a filtration on a
simplicial complex.  In other words, each $K_\varepsilon$ is a simplicial complex with
$K_\varepsilon\subseteq K_\eta$ for every $\varepsilon<\eta$, and we have
$K = \bigcup_{\varepsilon>0} K_\varepsilon$.  The $k$-th persistent homology of $ K $ is
given by
\[ \text{PH}_{k}(K):= \{ H_k(K_{\varepsilon}) \}_{\varepsilon\in \mathbb{R}_+} \] together
with the collection of linear maps
$ \psi^k_{\varepsilon,\eta}\colon H_k(K_{\varepsilon})\rightarrow H_k(K_{\eta}) $ induced by
the inclusion maps of $ K_{\varepsilon} \hookrightarrow K_{\eta} $ for all $k\in\mathbb{N}$
and $ \varepsilon<\eta$ in $\mathbb{R}_+ $.

\subsection{Bar codes}

Persistent homology produces a collection of intervals depending on the parameter
$ \varepsilon $ where we store the \emph{life-time} of topological features of the point
cloud via persistent homology.  Here by \emph{life-time} we mean the interval on which a
homology cycle is non-trivial as $\varepsilon$ ranges from 0 to $\infty$.  We record both
the \emph{birth}, i.e when a topological feature appears, and the \emph{death}, i.e when a
topological feature disappears, as $ \varepsilon $ increases. To illustrate the life-time,
we use \emph{barcodes} as introduced by~\cite{carlsson2005persistencebarcode}
and~\cite{ghrist2008barcodes}.

\begin{figure}[ht]
	\centering \includegraphics[scale=0.15]{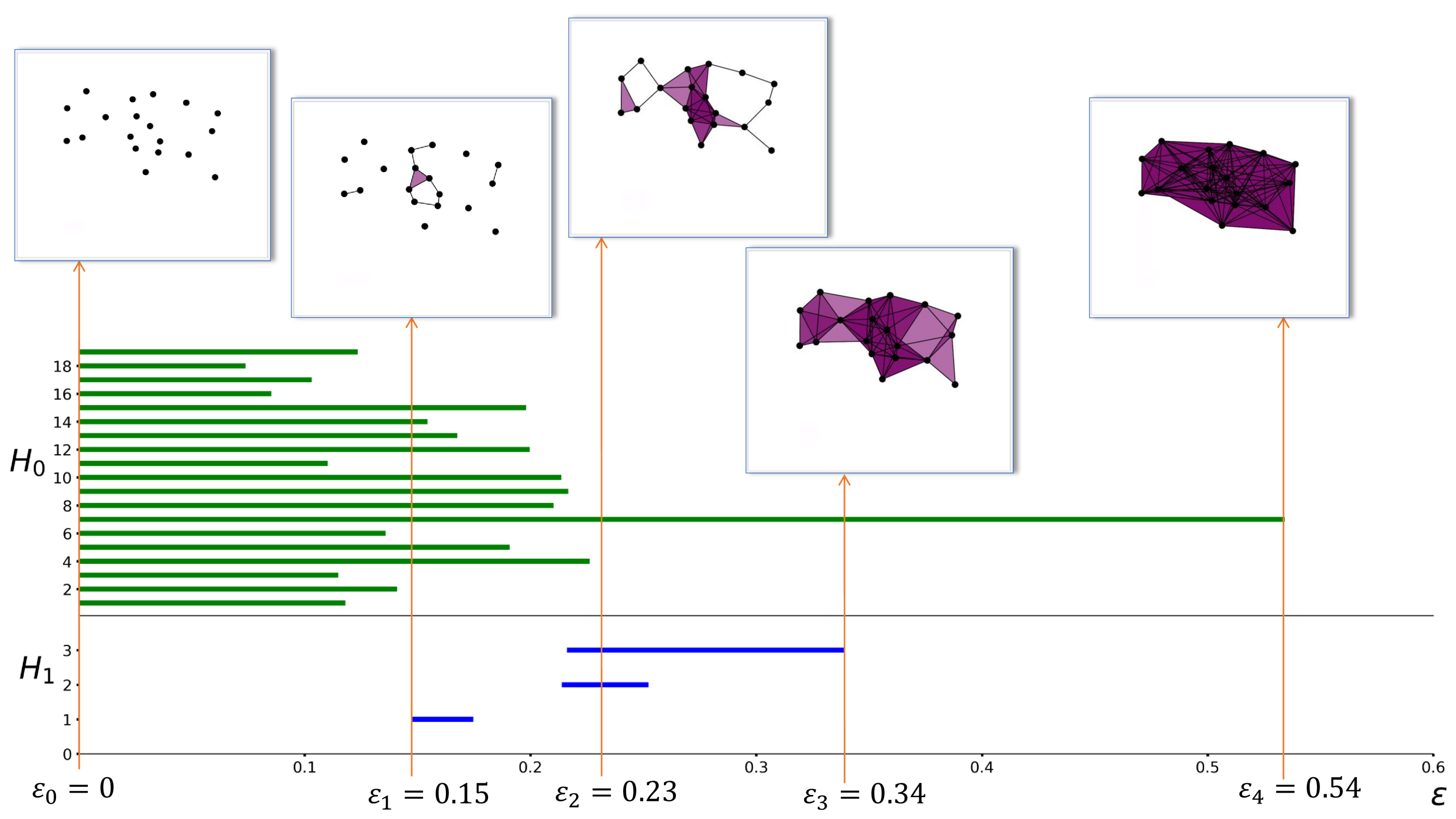}
	\caption{ An example barcode. } 
	\label{barcodeswithcomplex}
\end{figure}

In a barcode, we place the basis vectors for the homology on the vertical axis whereas the
horizontal axis represents the life span of each basis element in terms of the scale
parameter $ \varepsilon $. When we draw the vertical line at a particular $ \varepsilon_i$,
the number of intersecting line segments in a barcode is the dimension of the corresponding
homology group, i.e. the Betti number, for that parameter $\varepsilon_i$. See
Figure~\ref{barcodeswithcomplex}.

In Figure~\ref{barcodeswithcomplex}, one can see barcodes for zeroth and first persistent
homology together with the Vietoris-Rips complex corresponding to a particular
$ \varepsilon $.
For example, the blue horizontal line whose left endpoint on $ 0.22 $ and right endpoint is
on $ 0.25 $ represents a nonzero element in $ H_1(R_{0.22}) $ that persisted until
$ H_1(R_{0.25}) $ at which point it either disappeared or merged with another class.

We will postulate that the longest living topological features in the barcode are the
genuine topological features of the point cloud, whereas the shorter ones can be seen as
artificial artifacts of the method we use.  Notice also that there will always be one
connected component as $ \varepsilon $ grows large, i.e. the zeroth Betti number
$ \beta_{0} $ is always going to be 1 eventually.

\section{Hierarchical Clustering}\label{sect:HC}

Assume we have a connected metric space $(X,d)$, and let $\pi_0(X)$ be the set of connected
components of $X$. Assume we have a finite random sample of points $D\subseteq X$ taken from
$X$ whose distribution we do not know.  Our aim is to deduce any information about the set
of connected components of $X$ using $D$.  We are going to do this by finding a \emph{finite
	clustering} of $D$ which is a set function $c\colon D\to \mathbb{N}$ such that each
cluster $c^{-1}(i)$ lies within a distinct connected component for each $i\in\mathbb{N}$.

\subsection{Hierarchical clustering}

In its simplest form, in hierarchical clustering we have a function
$c_\varepsilon\colon D\to \mathbb{N}$ for each scale parameter $\varepsilon>0$.  This
function satisfies $c_\varepsilon(x)=c_\varepsilon(y)$ for any two points $x,y\in D$ when
there is a sequence of points $x_0,\ldots,x_m\in D$ such that $d(x_i,x_{i+1})< \varepsilon$
for every $i=0,\ldots,m-1$ where $x_0=x$ and $x_m=y$.  Notice that the clustering algorithm
is monotone in the sense that if $c_\varepsilon(x)=c_\varepsilon(y)$ then
$c_\eta(x)=c_\eta(y)$ for every $\eta>\varepsilon$.  Moreover, since $D\subseteq X$ is
finite and $X$ is connected, there is a large enough scale parameter $\varepsilon>0$ such
that the image of $c_\varepsilon$ is a single cluster.
\begin{figure}[ht]
	\centering
	\begin{algorithmic}
		\Procedure{Cluster}{$Data,\varepsilon$}
		\State $\mathcal{C}\gets \emptyset$
		\For{$x$ in $Data$}
		\State Add $\{x\}$ as a cluster to $\mathcal{C}$
		\EndFor
		\Repeat
		\State Find a distinct pair $(C_i,C_j)$ in $\mathcal{C}$ such that $d(C_i,C_j)<\varepsilon$
		\State Remove the clusters $C_i$ and $C_j$ from $\mathcal{C}$
		\State Add the new cluster $C_i\cup C_j$ to $\mathcal{C}$
		\Until{$d(C_i,C_j)\geq\varepsilon$ for all clusters}
		\State\Return $\mathcal{C}$
		\EndProcedure
	\end{algorithmic}
	\caption{Clustering function pseudocode.}\label{Algorithm1}
\end{figure}

\subsection{Linkage in hierarchical clustering}\label{sect:linkages}

As we increase the scale parameter $\varepsilon>0$ we start forming \emph{clusters} of
points.  Since we replace points with clusters, we are going to need to calculate
distances between clusters.  See~Algorithm~\ref{Algorithm1}.

For a fixed $\varepsilon>0$, let us use $C_i=c^{-1}_\varepsilon(i)$ to denote a cluster, and
set $n_i = |C_i|$.  Let us use $d_{ij}$ for the distance between the cluster $ C_i $ and
$ C_j $. Lance and Williams \cite{lance1967general} introduced to the following general
formula for calculating distances between clusters
\[ d_{(ij)k} = \alpha_{ijk} d_{ik} +\alpha_{jik} d_{jk}+ \beta_{ijk} d_{ij} +
  \gamma|d_{ik}-d_{jk}| \] for parameters $\alpha_{ijk}$, $\beta_{ijk}$ and $\gamma$ to be
determined.  Here, $ d_{(ij)k} $ denotes the distance between the clusters $ C_{k} $ and
$ C_{ij} = C_i \cup C_j $ which is merged in a single cluster.  We list the parameters for
commonly used methods of calculating distances between clusters in
Table~\ref{table:distance_methods}. See~\cite{lance1967general} for details.

\begin{table}[ht]
	\centering
	\caption{Commonly used methods to determine $d_{(ij)k}$.}\label{table:distance_methods}
	\setlength{\arrayrulewidth}{0.25mm}
	\setlength{\tabcolsep}{10pt}
	\setlength{\extrarowheight}{6pt}
	\renewcommand{\arraystretch}{1.5}
	\begin{tabular}{lccc}
		Linkage  & $\alpha_{ijk} $ & $\beta_{ijk}$ & $\gamma$\\ \hline
		Single   & $\frac{1}{2}$ & 0 & $-\frac{1}{2}$ \\ 
		Complete & $\frac{1}{2}$ & 0 & $\frac{1}{2}$ \\ 
		Average  & $\displaystyle\frac{n_i}{n_i + n_j}$ & 0 & 0 \\ 
		Ward     & $\displaystyle\frac{n_i + n_k}{n_i + n_j + n_k}$ & $\displaystyle\frac{-n_k}{n_i + n_j + n_k}$ & 0 \\ \hline
	\end{tabular}
\end{table}

  \subsection{Evaluation of clusters}\label{sect:silhouette}

  Performance evaluations in supervised learning tasks are easier and more understandable
  than unsupervised learning tasks such as clustering. However, there are also some useful
  metrics one can use to evaluate clustering tasks. In this subsection, we are going to
  review the metrics we are going to use to evaluate our classification and clustering
  algorithms.  Particular metrics we are going to review are the mutual
  information~\cite{mutual}, Rand index~\cite{randscore}, and homogeneity,
  completeness~\cite{homegcomplete}, and the silhouette scores~\cite{silhouette}.

  For this subsection, assume $X$ is the random variable that designates the true labels
  $\mathcal{X} = \{x_1,\cdots,x_n\}$ while $Y$ is the random variable that designates the
  predicted labels $\mathcal{Y} = \{y_1,\cdots,y_m\}$ of a classification task. Let
  $\tau(d)$ be the true label and $\pi(d)$ be the predicted label of a data point $d$. We
  use $H(Z)$ to denote the Shannon entropy of a discrete random variable $Z$ which is
  defined as $$H(Z) = -\sum_{z\in Z} p(z)\log p(z).$$

  \paragraph{Mutual information:}

  The mutual information of the pair of the random variables $ X $ and $ Y $ is given as
$$ MI(X,Y) = H(X) + H(Y) - H(X,Y) $$
Let $\mathrm{C}_x=\tau^{-1}(x)$ be the set of samples with true label $x\in\mathcal{X}$ of
size $n_x$, and let $\mathrm{C}_y=\pi^{-1}(y)$ be the set of samples with predicted label
$y\in\mathcal{Y}$ of size $n_y$. The joint distribution $ p(X=x,Y=y) $ which is an
observation drawn at random falls into clusters $ \mathrm{C}_x $ and $ \mathrm{C}_y $ turn
out to be $ \frac{|\mathrm{C}_x \cap \mathrm{C}_y |}{n} $ with the marginal probability
$ p(X=x) = \frac{n_x}{n} $ and $ p(Y=y)=\frac{n_y}{n} $. So, the mutual information of the
pair $(X,Y)$ can be written in terms of these cardinalities as
$$ MI(X,Y) = \sum_{x\in\mathcal{X}} \sum_{y\in \mathcal{Y}} \frac{|\mathrm{C}_x \cap \mathrm{C}_y |}{n} \log \left( \frac{n |\mathrm{C}_x \cap \mathrm{C}_y |}{n_x n_y}  \right). $$

\paragraph{Homogeneity and completeness:}

If all clusters contain only data points that are members of a single class, then we say
that our clusters are homogeneous. On the other hand, if the data points that are members of
a given class are elements of the same cluster then we say that our clusters are complete.
Formally, the homogeneity and completeness scores are respectively defined by:
$$
hom(X,Y) = 1-\frac{H(X\mid Y)}{H(X)} \quad\mbox{ and }\quad
comp(X,Y) = 1-\frac{H(Y \mid X)}{H(Y)}
$$
Notice that both homogeneity and completeness scores are not affected by permutations of
labels.

\paragraph{Rand index:}

The Rand index counts sample of unordered distinct pairs of data points $\{u,v\}$ for which
$\tau$ and $\pi$ agree, and also disagree. Formally, the Rand index is defined as
$$ \mathrm{RI}=\frac{a+b}{\binom{n}{2}}, $$
where 
$$ a = | \{\{u,v\}: \tau(u)=\tau(v)\text{ and } \pi(u)=\pi(v) \} | $$
and
$$ b = | \{\{u,v\}: \tau(u)\neq\tau(v)\text{ and } \pi(u)\neq\pi(v) \} | $$
and $\binom{n}{2}$ is the total number of different possible unordered pairs in the dataset.

\paragraph{Silhouette score:}

We are going to use the silhouette score as defined in~\cite{silhouette} to evaluate a
clustering model on a dataset.

Assume we have a set of clusters $C_1,\ldots,C_k$, and a dissimilarity measure $d(x,y)$ for
every pair of points in our dataset. Let $U(x) = \pi^{-1}\pi(x)$ be the cluster that $x$
belongs to and let
\[ d(x,C_j) = \frac{1}{|C_j|} \sum_{y\in C_j} d(x,y) \] The silhouette score $s(x)$ of a
point $x$ in our dataset is given by the ratio
\[ s(x) = \frac{a(x) - b(x)}{\max(a(x),b(x))} \] where
\[ a(x) = d(x,U(x))\text{ and } b(x) = \min_{C\neq U(x)} d(x,C). \] Silhouette scores take
values in the interval $[-1,1]$ and values closer to 1 indicate that clusters are
\emph{well-formed}: they show low average intra-class similarity while maintaining a high
average inter-class dissimilarity.

\section{A bridge between persistent homology and hierarchical clustering}\label{sect:PHandHC}

\subsection{Cophenetic matrix}

An important notion we need in studying and comparing clustering methods is \emph{the cophenetic matrix} as defined in~\cite{sokal1962comparison,cophenetic1973numerical,JainDubes88}.

Assume we have a clustering function $c_\varepsilon \colon D\to \mathbb{N}$, and let
$\mathcal{C} = \{C_i = c_\varepsilon^{-1}(i)\mid i\in \mathbb{N}\}$.  Let
$ \varepsilon_{ij} $ be the proximity level at which the clusters $C_i$ and $C_j$ merge to
form $ C_{ij} $ for the first time. We record these numbers in the cophenetic matrix
$C_\varepsilon(D) = (\varepsilon_{ij})$ for any pair of clusters $C_i$ and $C_j$.  The
cophenetic distance is a metric under the assumption of monotonicity~\cite{PatternRecognitionCophenetic}.

\subsection{Homological cophenetic distance}\label{sect:homological-distance}

Given a point cloud $D$, we consider the Vietoris-Rips complex $R_\varepsilon(D)$.  By
gradually increasing $ \varepsilon $, we get a filtered simplicial complex, and thus, we can
calculate the persistent homology associated with this filtration.

Recall that when we have a filtered simplicial complex $\{R_\varepsilon\}_{\varepsilon>0}$,
we have homology groups $\{H_k(R_\epsilon)\}_{\varepsilon}$ and connecting linear maps
$\psi^k_{\varepsilon,\eta}\colon H_k(R_\epsilon)\to H_k(R_\eta)$ for every pair
$\varepsilon<\eta$ and for every $k\in\mathbb{N}$.  We would like to emphasize that even
though $R_{\varepsilon}\subseteq R_{\eta}$ the induced maps in homology
$\psi^k_{\varepsilon,\eta}$ need not be injections.

Now, for each linearly independent pair of homology classes $ \alpha$ and $ \beta $ in
$ H_{k}(R_{\varepsilon}) $ one can test if $\psi^k_{\varepsilon,\eta}(\alpha)$ and
$\psi^k_{\varepsilon,\eta}(\beta)$ are still linearly independent in $ H_{k}(R_{\eta}) $.
If the pair $ \psi^n_{\varepsilon,\eta}(\alpha)$ and $ \psi^n_{\varepsilon,\eta}(\beta)$ is
non-zero and fails to be linearly independent then we will say that two classes $ \alpha $
and $ \beta $ merged at time $ \eta $.  

\begin{definition}\label{defn:1}
	\emph{$k$-th homological cophenetic distance} is defined as
	\begin{equation*}
		\label{eq:1}
		D_k(\alpha,\beta) =
		\inf\left\{\eta-\varepsilon\geq 0 \mid \psi^k_{\varepsilon,\eta}(\alpha) \text{ and }\psi^k_{\varepsilon,\eta}(\beta ) \text{ are non-zero and linearly dependent} \right\}. 
	\end{equation*}
\end{definition}
\begin{proposition}
	The cophenetic metric is non-archimedean:
	\begin{equation*}
		\label{prop:1}
		D_k(\alpha,\beta)\leq \max(D_k(\alpha,\gamma),D_k(\gamma,\beta))
	\end{equation*}
	for every $\alpha,\beta,\gamma\in H_k(R_\varepsilon)$.
\end{proposition}

\subsection{Non-archimedean metrics and hierarchical clustering}

It has been known that hierarchical clustering methods and non-archimedean metrics are
intimately related~\cite{JardineSibson71,Hartigan85,JainDubes88}. But it was Carlsson
and Memoli who proved that hierarchical clustering methods and non-archimedean metrics
are naturally equivalent in~\cite{carlsson2010hierarchical}.  The immediate corollary
is that the homological cophenetic metric we defined in Definition~\ref{defn:1} does
correspond to a unique hierarchical clustering scheme on every homology group.

Recall that the connected components of a topological or a metric space $X$ are encoded
in the zeroth homology classes of $X$. This means, one can naturally compare various
clustering schemes on a dataset embedded in $X$ with the homological cophenetic
distance for the zeroth homology by varying the underlying metric of $X$ if the ambient
space $X$ allows it.  The rest of the paper will follow this route.

Observe also that we have homological cophenetic distance for all homology groups, not
just zeroth homology, which do not have obvious connections with the point clouds that
form a dataset.  Investigating the ramifications of the non-archimedean homological
cophenetic distance for higher homology groups is going to be the subject matter of a
subsequent paper.

\subsection{The zeroth homology and hierarchical clustering}

For the zeroth homology there is a simplification: all homology classes appear at
$\varepsilon=0$ since every point is a connected component by itself, and they never
disappear only merge as $\varepsilon$ goes to $\infty$ for the subsequent Vietoris-Rips
complexes.  Thus, it is enough to test whether the classes $\psi^0_{0,\varepsilon}(\alpha)$
and $\psi^0_{0,\varepsilon}(\beta )$ are linearly independent in $H_{0}(R_{\varepsilon})$ as
$\varepsilon$ varies.

We also note that each point $x\in D$ is a homology class in $H_0(R_0)$, and these points
mark the rows and columns of the cophenetic matrix $C_0(D)$ coming from hierarchical
clustering.  As a consequence we can compare these matrices.

\subsection{Mantel Test}\label{mantelsection}

As we stated above, we need to compare different dendrograms, or equivalently, different
distance matrices on a given dataset.  For this purpose we are going to use the Mantel test as defined in~\cite{mantel1967detection}. It is commonly used in biology and ecology to
compare phylogenetic trees. Mantel test is a non-parametric statistical method that computes the
significance of correlation between rows and columns of a matrix through permutations
of these rows and columns in one of the input distance matrices.

We consider two distance or cophenetic matrices $D_1=(x_{ij})$ and $D_2=(y_{ij})$ of size
$n\times n$. The normalized Mantel statistic $ r $ is defined as
\[ r = \frac{2}{(n-2)(n+1)}\sum_{i=1}^{n}\sum_{j=i+1}^{n}
\left(\frac{x_{ij}-\bar{x}}{s_x}\right) \left( \frac{y_{ij}-\bar{y}}{s_y}\right) \]
where
\begin{enumerate}[(i)]
	\item $ \bar{x} $ and $ \bar{y} $ are averages of all entries of each matrix, and
	\item $ s_x $ and $ s_y $ are the standard deviations for $ x $ and $ y $.
\end{enumerate}
The test statistic is the Pearson product-moment correlation coefficient $ r \in[-1,1] $.
Having a value in the neighborhood of $ -1 $ indicates strong negative correlation whereas
$ +1 $ indicates strong positive correlation, and 0 indicates no relation.

In order to estimate the sampling distribution of the standardized Mantel statistic under
the null-hypothesis (no correlation between the distance matrices), random permutations of
the rows (or equivalently columns) of the distance matrices are used to get a set of values
of the statistic.  Then whether the null-hypothesis is rejected depends on the value of the
Mantel statistic: If the calculated statistic is unlikely to have been obtained under the
null-hypothesis then the null-hypothesis is rejected.  See~\cite[Sect. 10.5]{Legendre2012}
for details.

\section{Experiments}\label{sect:experiments}

To determine if our research is sound, we performed numerical experiments\footnote{The
  computational tools we use in this section are as follows: To compare dendrograms, we use
  tools come from the \textbf{dendextend}~\cite{galili2015dendextend} and
  \textbf{vegan}~\cite{veganRpackage} packages of the R programming language~\cite{R}. For
  the map of Turkey, we used Generic Mapping Tools~\cite{gmt6visualization}. To compute
  cophenetic distance matrix we used SageMath~\cite{sagemath}. In order to compute and
  visualize clusters, we used python programming language~\cite{python} and its
  \textbf{scikit-learn} library~\cite{sklearn}.}  on synthetic and different real
datasets. Statistically, we compare the dendrograms we obtained from the Euclidean distance
and the dendrograms we obtained from the cophenetic distance for the zeroth persistent
homology. Also, this statistically comparison methods is done for various distance metrics,
too.  On the other hand, we evaluate the clustering measurements from hierarchical
clustering algorithm by changing the metrics on different datasets. In this section, we are
going to summarize these experiments.

\subsection{A sample of cities in Turkey}\label{subsect:TurkishCities}

\begin{figure}[ht]
	\centering
	\includegraphics[scale= 0.85]{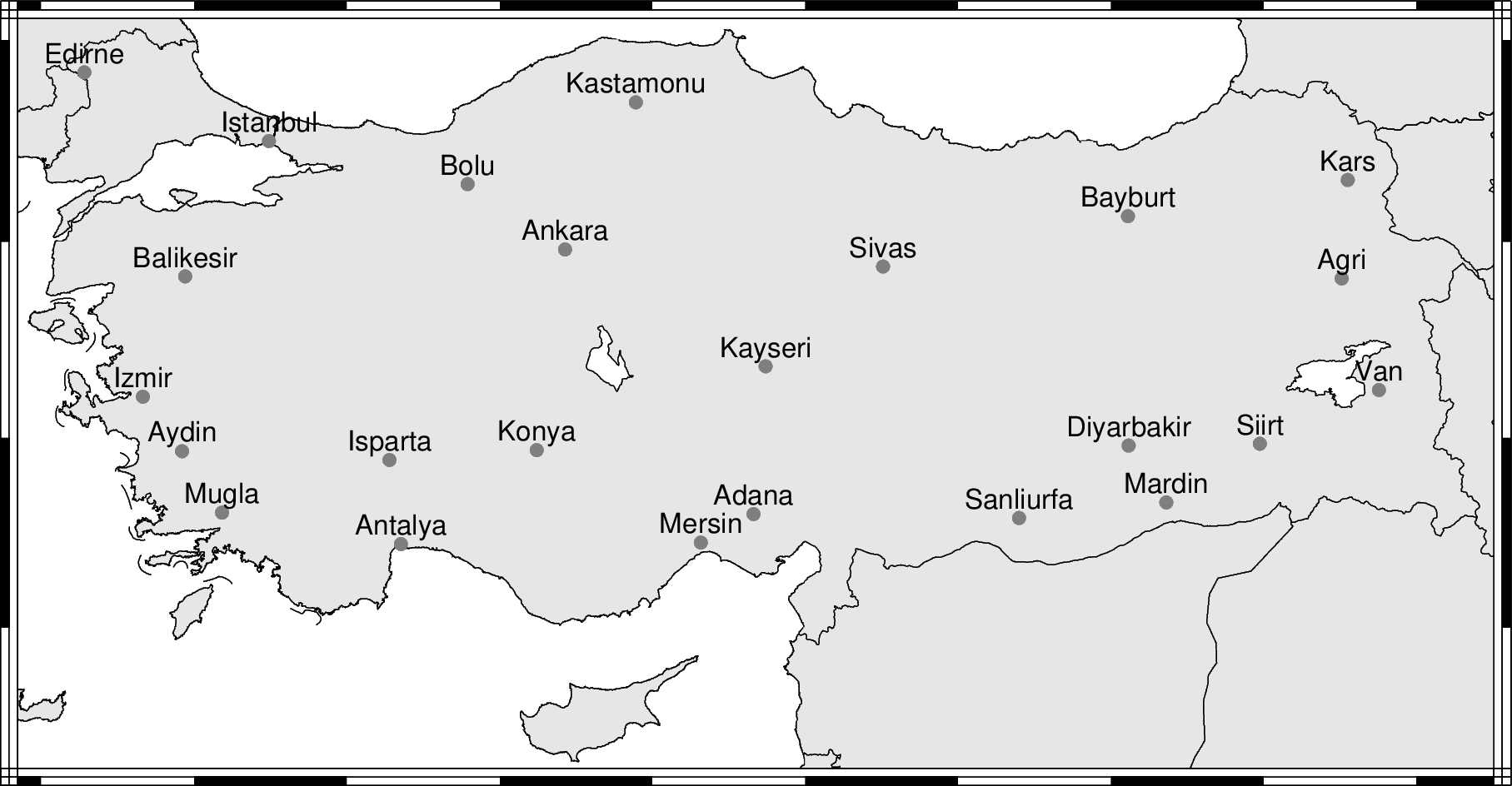}
	\caption{A sample of cities in Turkey}
	\label{turkish_cities}
\end{figure}

For our first experiment, we used a subset 24 of cities in Turkey whose coordinates are
encoded as longitudes and latitudes in radians. See Figure~\ref{turkish_cities}.

\subsubsection{Bar codes and dendrograms}

\begin{figure}[ht]
	\centering
	\includegraphics[scale=0.15]{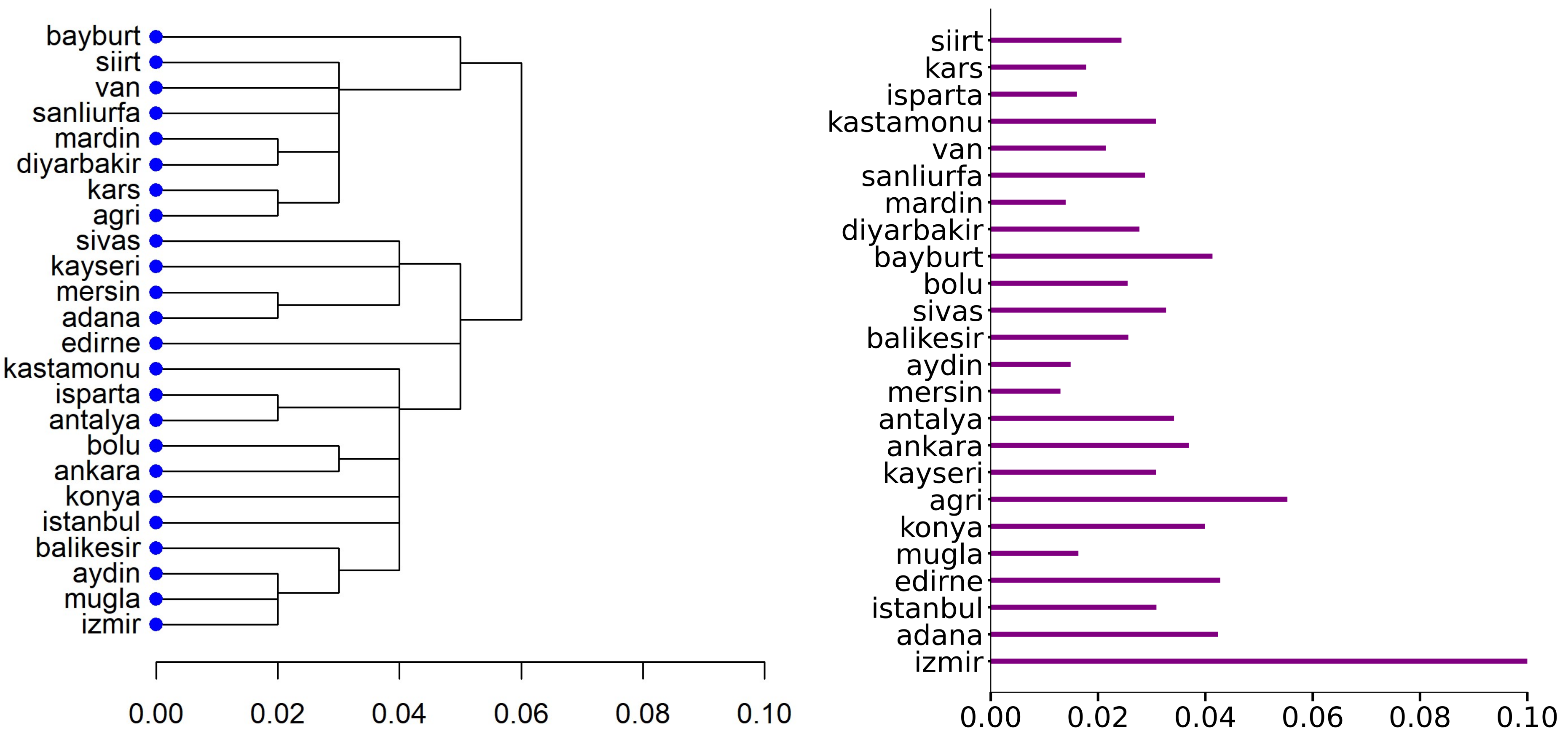}
	\caption{Hierarchical enriched barcodes and classical barcodes in TDA }
	\label{twobarcodestogether}
\end{figure}

The left hand side of Figure \ref{twobarcodestogether} is the dendrogram we obtained from
cophenetic homological distance matrix for the zeroth homology.  The right hand side of
Figure \ref{twobarcodestogether} is the ordinary barcode obtained from the zeroth persistent
homology which displays the birth and death times of each homology class, whereas the left
hand side is the dendrogram that indicates which classes merge.

\subsubsection{Comparison of Dendrograms}\label{subsect:dendrogram-comparison}

\begin{figure}[ht]
	\centering
	\includegraphics[width=\textwidth, height=7.5cm]{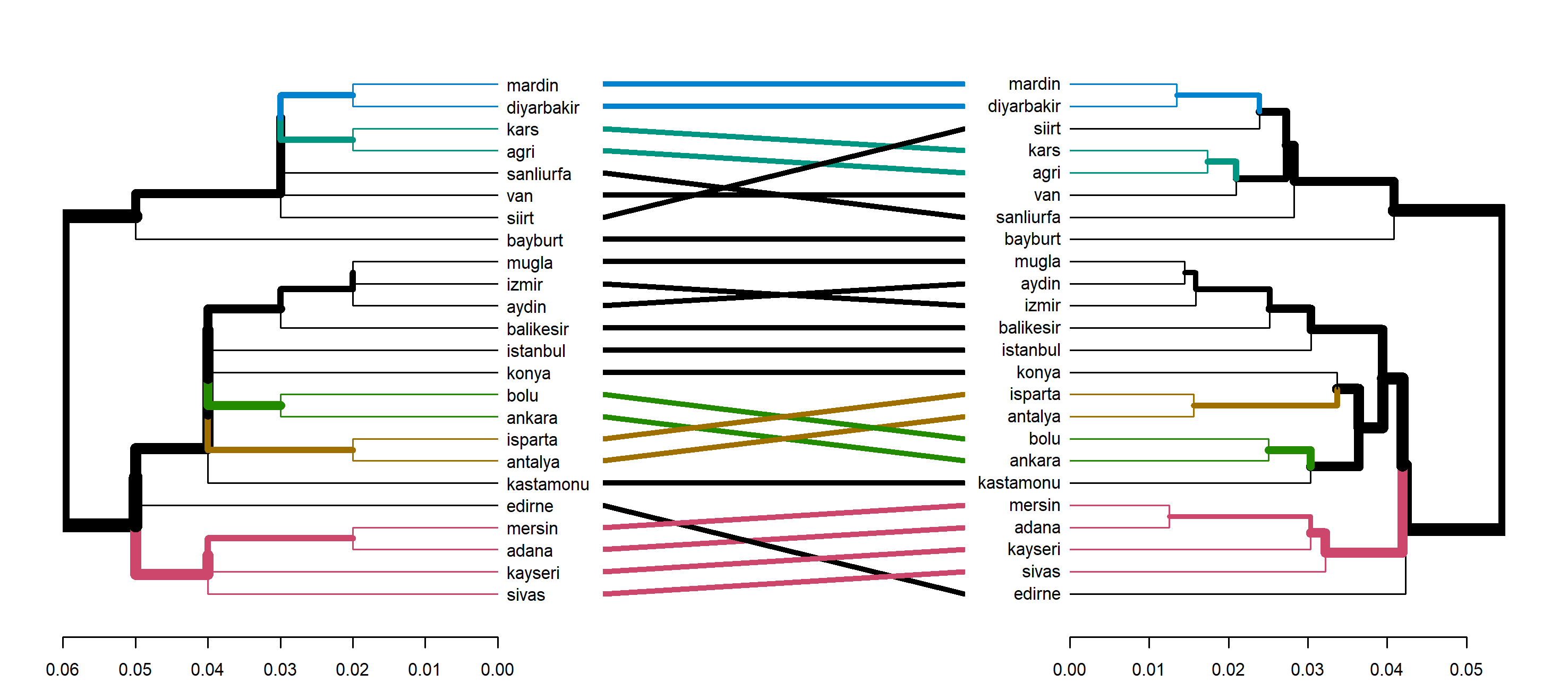}
	\caption{Tanglegram of the dendrograms of the cophenetic (left) and Euclidean
		distances (right).}
	\label{tanglegram}
\end{figure}

Next, we apply the hierarchical clustering (with single linkage), using the Euclidean
distance matrix $ E(D) $, and the homological cophenetic distance matrix $ C_0(D) $ for
the zeroth persistent homology. The resulting dendrograms are given in
Figure~\ref{tanglegram}. We also align the labels from both dendrograms without changing the
underlying cluster structure. In a tanglegram representation, one compares the tree
structures using a metric derived from matches between labels placed on
branches~\cite{scornavacca2011tanglegrams,fernau2010comparing,buchin2008drawing}.

For the next phase, we need to compare dendrograms. We are going to use the Mantel test (See
Section~\ref{mantelsection}) for this task.  The resulting statistic is a measure of how
well the labels of the two dendrograms are aligned.  For the sample of cities we used, the
Mantel statistic value we obtained for the matrices $E(D)$ and $C_0(D)$ was $ 0.98 $ with
p-value of $0.001$.

\subsubsection{A full comparison of metrics}

In Section~\ref{subsect:dendrogram-comparison}, we explicitly compared cophenetic
homological metric with the Euclidean metric.  In this section, we extend the
comparison to a variety of metrics.  As we discussed in Section~\ref{mantelsection}, we
are going to use the Mantel test for this comparison.

For our comparisons we used the $L^2$-metric (also known as the \emph{Euclidean metric}),
the $L^1$-metric (also known as \emph{taxicab distance}, \emph{city-block distance}, or the
\emph{Manhattan distance}), $L^p$-metric (also known as the \emph{Minkowski distance}) with
$p=1/2$, the cosine similarity converted to a dissimilarity function, and Bray-Curtis
dissimilarity~\cite{braycurtis1957}, and \cite[Eq. 7.58]{Legendre2012}. Our resulting Mantel
test table is given in Table~\ref{table:mantel}.

\newlength{\width}
\width4mm

\begin{table}[ht]
	\caption{Mantel statistics of metrics on the cities of Turkey dataset.}\label{table:mantel}
	\centering
	\setlength{\tabcolsep}{2pt}
	\footnotesize
	\begin{tabular}{ p{5\width} p{\width} p{\width} p{\width} p{\width} p{\width} p{\width}}
		\multicolumn{7}{c}{{PAIRWISE MANTEL STATISTICS}}\\\hline
		{\bf Metrics}
		&\multicolumn{1}{c}{{\bf Bray-Curtis}}
		&\multicolumn{1}{c}{{\bf Cosine}}
		&\multicolumn{1}{c}{{\bf Manhattan}}
		&\multicolumn{1}{c}{{\bf Euclidean}}
		&\multicolumn{1}{c}{{\bf Minkowski}}
		&\multicolumn{1}{c}{{\bf Homological}}\\\hline
		\bf Bray-Curtis & 1.00 & 0.64 & 0.96 & 0.90 & 0.90 & 0.90  \\\hline
		\bf Cosine      &      & 1.00 & 0.61 & 0.52 & 0.69 & 0.59   \\\hline
		\bf Manhattan   &      &      & 1.00 & 0.96 & 0.87 & 0.97  \\\hline
		\bf Euclidean  &      &      &      & 1.00 & 0.75 & 0.98 \\\hline
		\bf Minkowski   &      &      &      &      & 1.00 & 0.78  \\\hline
		\bf Homological &      &      &      &      &      & 1.00  \\\hline
	\end{tabular}\normalsize
\end{table}

The results indicate that our homological cophenetic distance, produce results most similar
to the Euclidean metric and the Manhattan metric, and is most dissimilar to the cosine
similarity. Homology is a topological invariant, and therefore, is
  impervious in any perturbations of the underlying metric structure. However, the
  filtration structure on the Rips complexes we used rely heavily on the underlying metric.
  Thus one can surmise that, the main contributing factor for this similarity might be the
fact that we used Euclidean metric when we formed Rips complexes. 

\subsection{Cophenetic distance on different datasets}\label{subsect:OtherDatasets}~        
\begin{table}[ht]
	\centering
	\caption{Datasets used and their properties.}
	\setlength{\tabcolsep}{4pt}
	\footnotesize
	\begin{tabular}{lcccc}
		\bf Dataset & \bf \#Instances & \bf \#Attributes & \bf Supervised & \bf \#Classes\\
		\hline Turkish Cities  & 82 & 2 & \text { No } & - \\
		\hline Iris & 150 & 4 & \text { Yes } & 3 \\
		\hline Cancer Coimbra & 116 & 10 & \text { Yes } & 2 \\
		\hline Synthetic (total separation) & 100 & 100 & \text { Yes } & 4 \\
		\hline Synthetic (with mixture) & 100 & 2 & \text { Yes } & 4 \\
		\hline
	\end{tabular}
	\normalsize
\end{table}

In the previous section, we compared the dendrograms coming from different metrics on the
coordinates of a small sample cities in Turkey.  In this section, we are going to evaluate
the clusters coming from hierarchical clustering algorithms by varying the metrics on 5
different datasets: all cities in Turkey, the \emph{Iris
  Dataset}~\cite{iris}, 
the \emph{Cancer Coimbra Dataset}~\cite{cancer} and two synthetic
  datasets that we generated. One of these synthetic datasets has 4 linearly separable
  clusters, and the other contains 4 clusters with mixing along their boundaries. For the
synthetic datasets we used {\tt make\_blobs} function from the scikit-learn
library~\cite{sklearn} of the python language~\cite{python}.

\subsubsection{Silhouette scores}

\begin{figure}[ht]
	\centering
	\begin{subfigure}[b]{.49\textwidth}
		\centering
		\includegraphics[width=\textwidth]{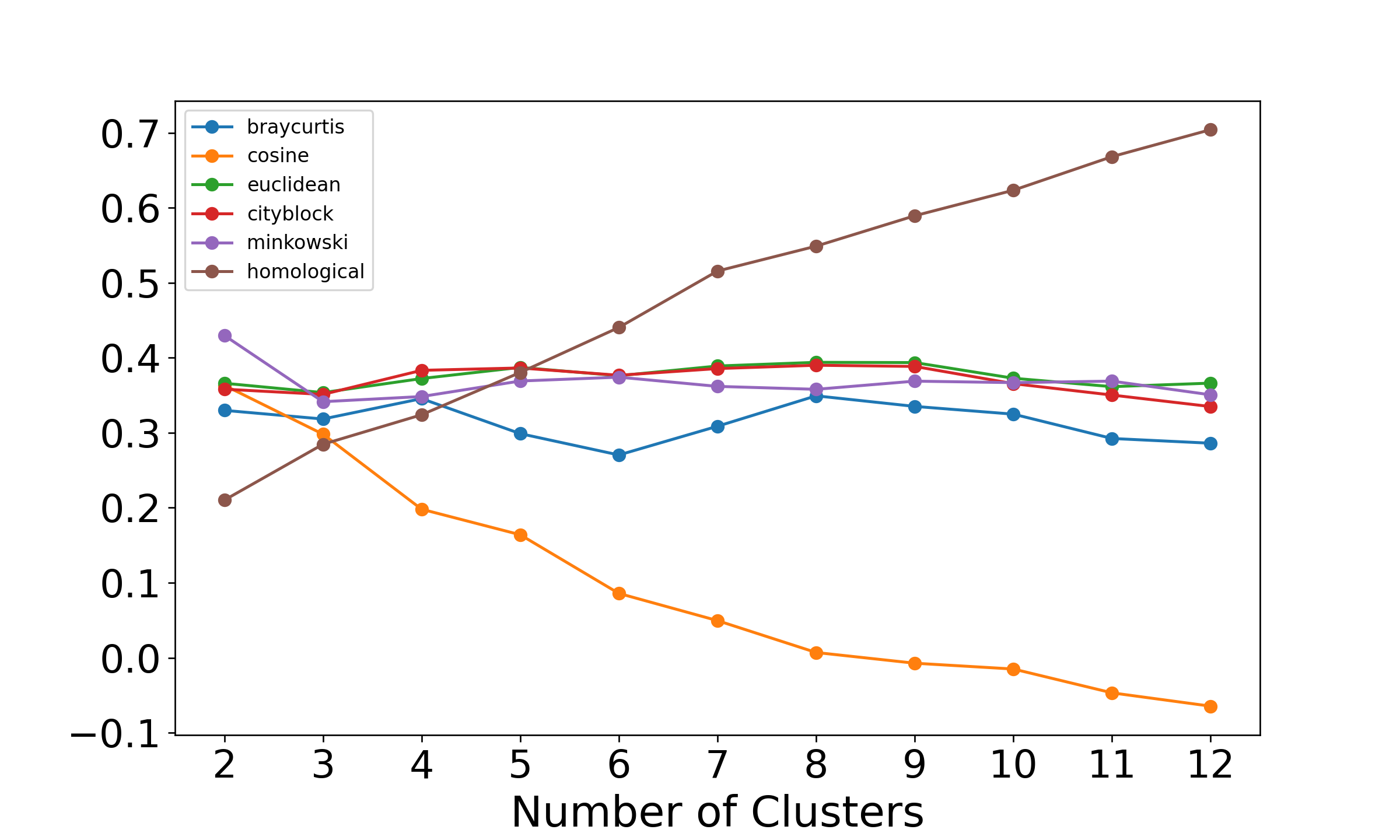}
		\caption{Turkish cities.}
	\end{subfigure}
	\hfill
	\begin{subfigure}[b]{.49\textwidth}
		\centering
		\includegraphics[width=\textwidth]{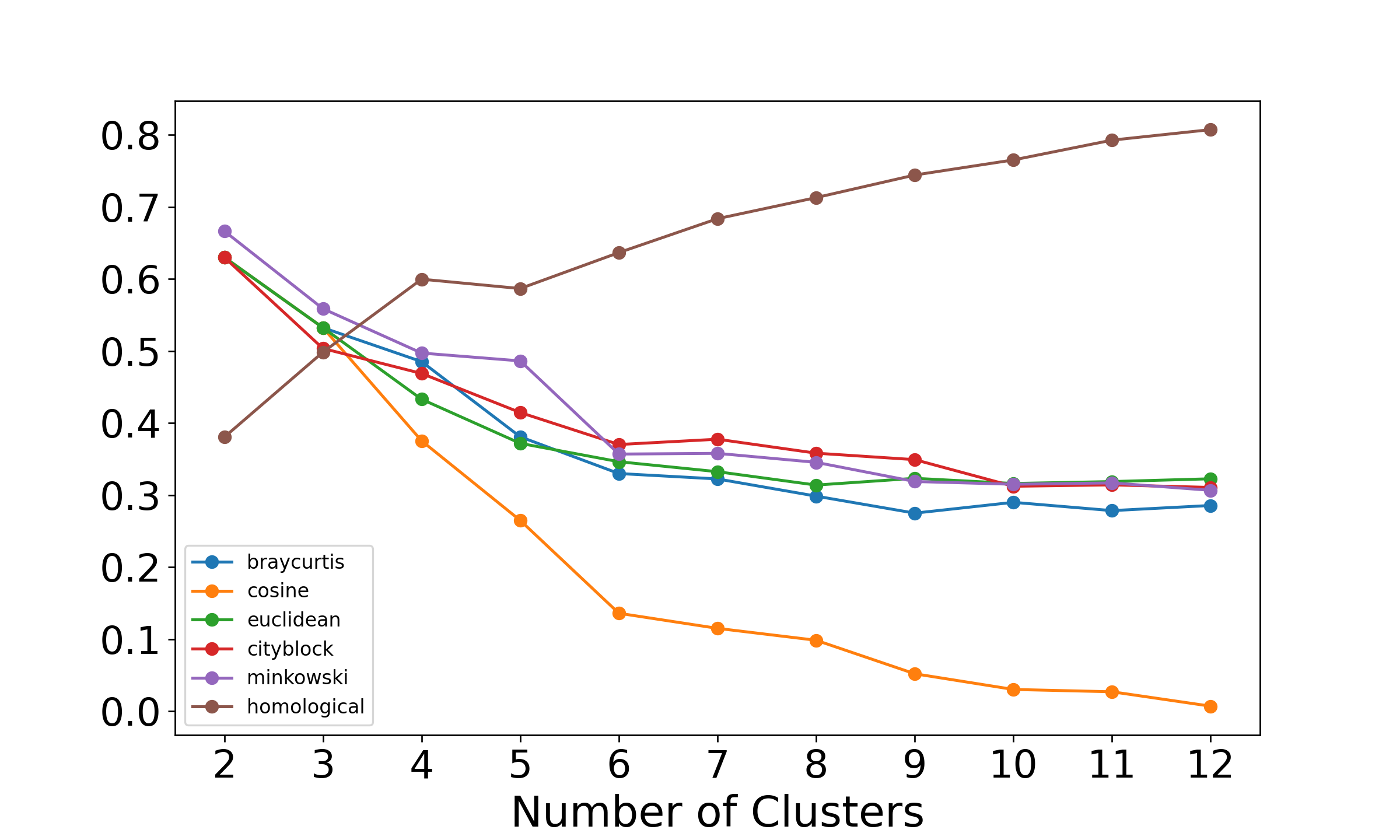}
		\caption{Iris dataset.}
	\end{subfigure}
	\hfill
	\newline
	\begin{subfigure}[b]{.49\textwidth}
		\centering
		\includegraphics[width=\textwidth]{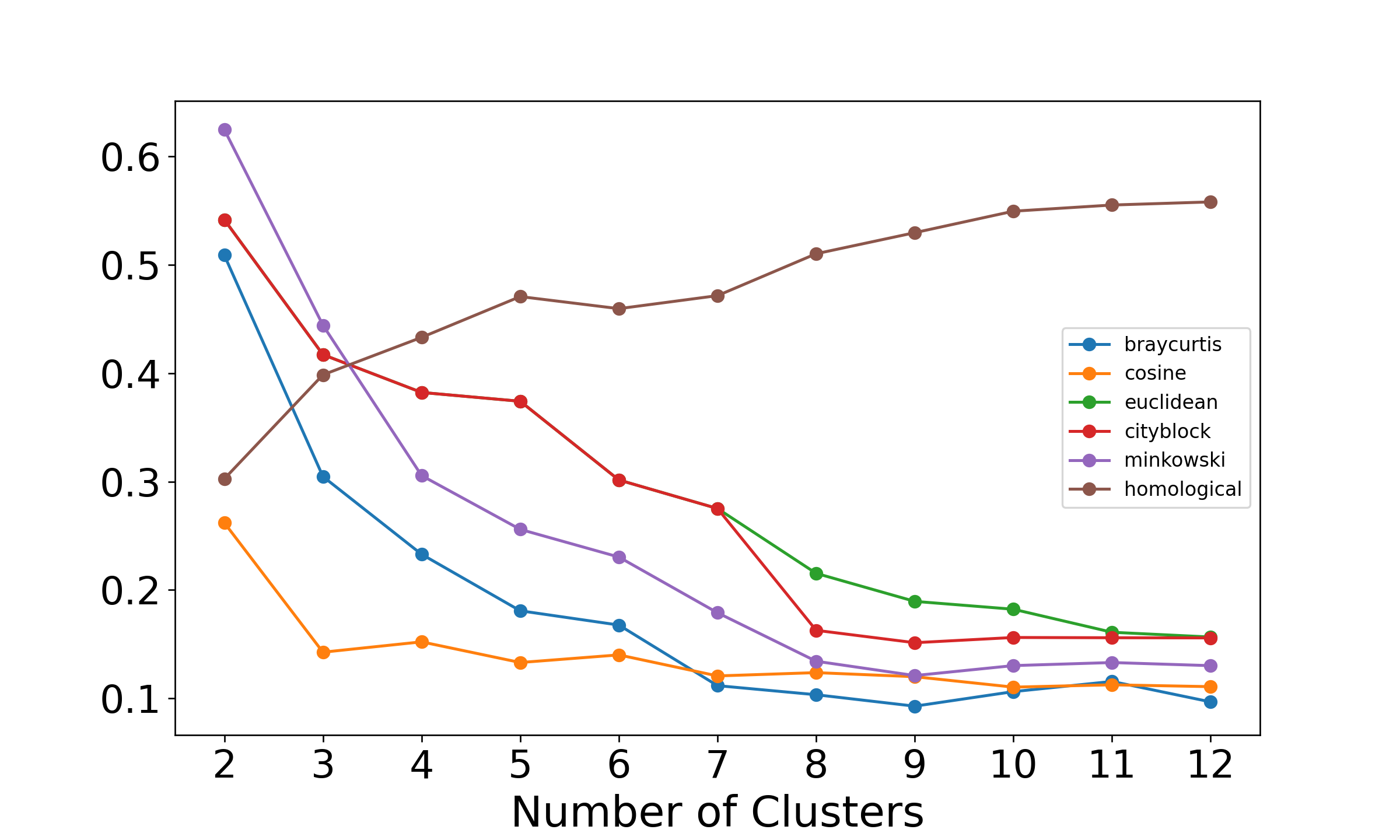}
		\caption{Cancer dataset.}
	\end{subfigure}
	\hfill
	\begin{subfigure}[b]{.49\textwidth}
		\centering
		\includegraphics[width=\textwidth]{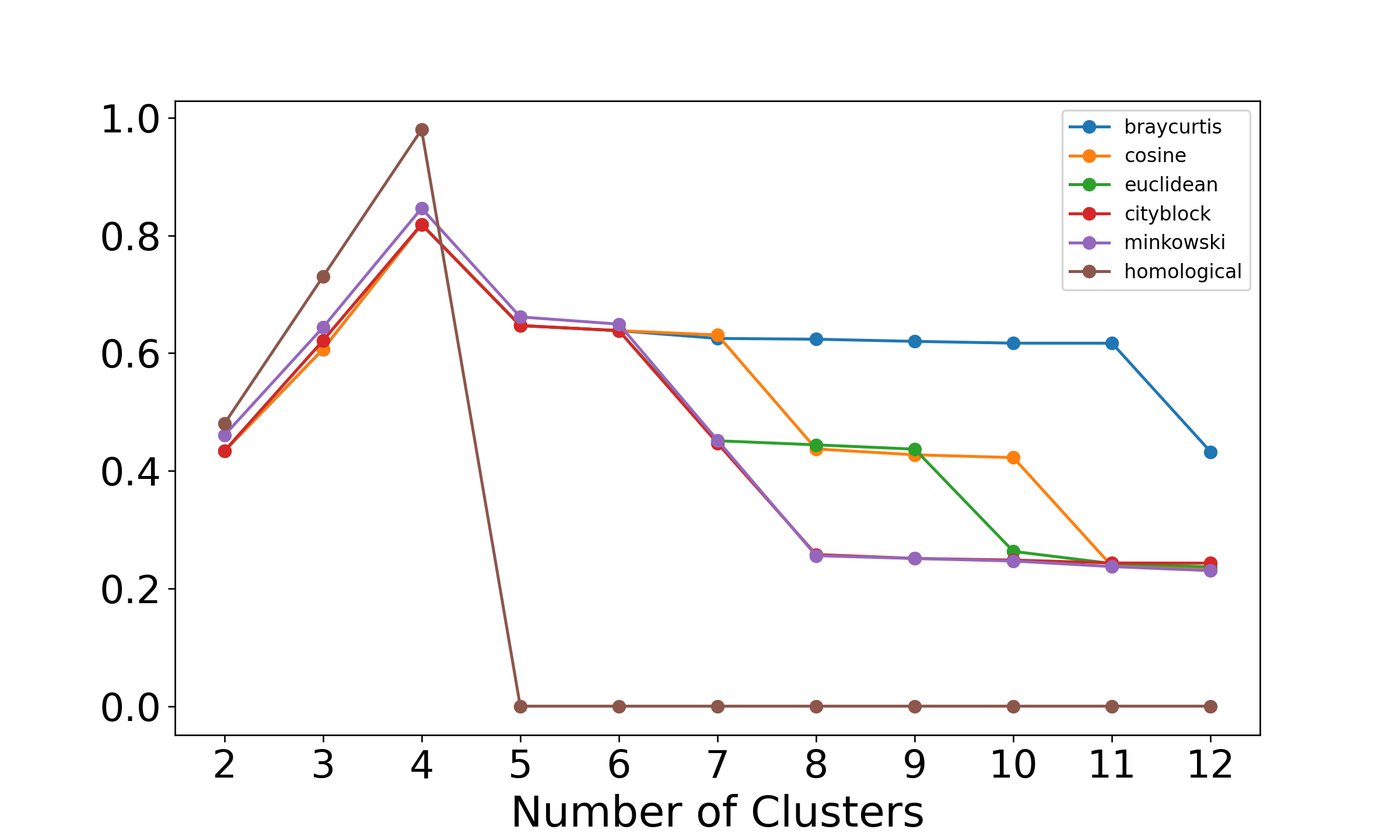}
		\caption{Separable synthetic dataset.}
	\end{subfigure}
	\hfill 
	\begin{subfigure}{.49\textwidth}
		\centering
		\includegraphics[width=\textwidth]{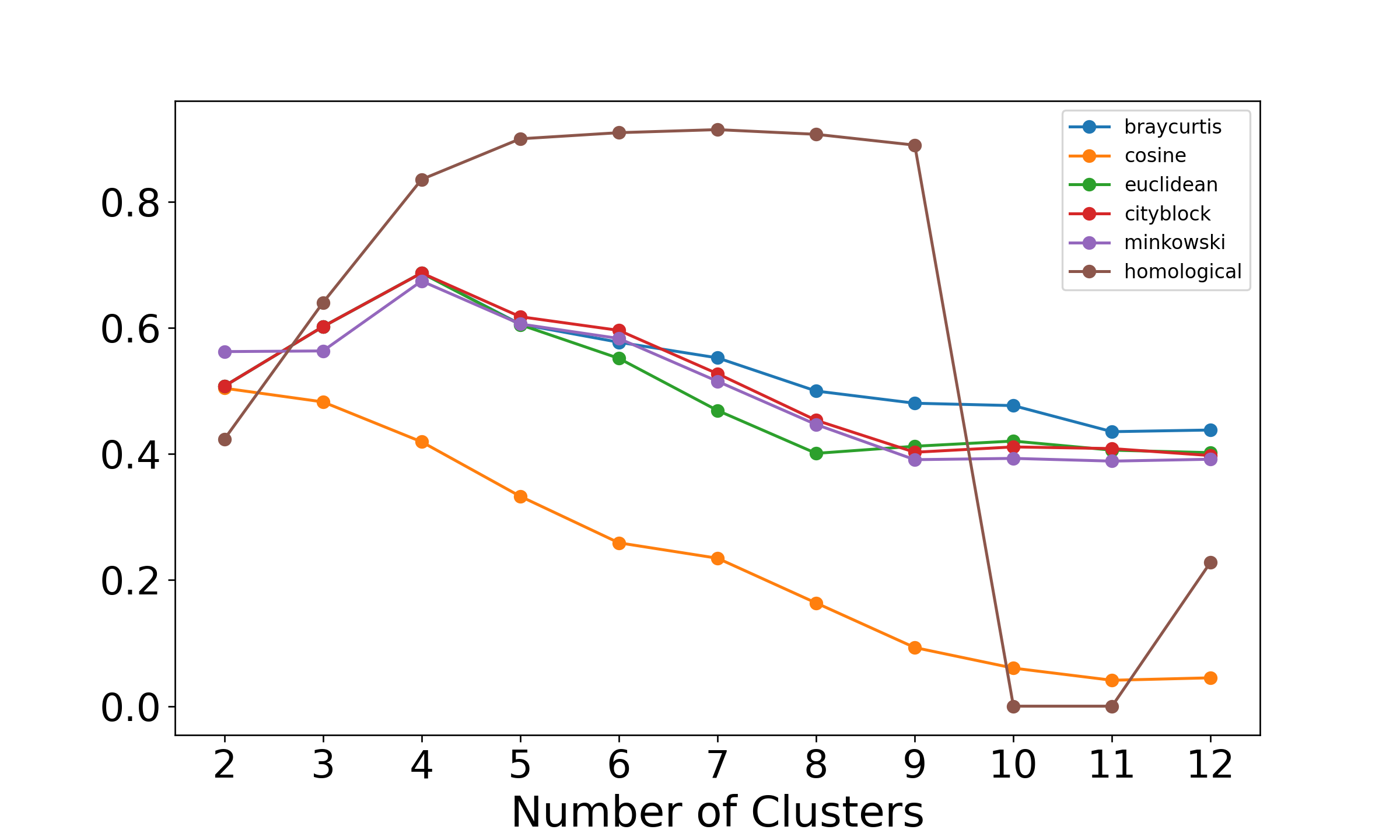}
		\caption{Synthetic dataset with mixing.}
	\end{subfigure}
	\hfill 
	\caption{Silhouette scores for each dataset.}
	\label{fig:silhouette graphs}
\end{figure}

In order to compare the resulting clusters on each dataset, we are going to use silhouette
scores as discussed in Section~\ref{sect:silhouette}.  The results of average silhouette
scores for different number of clusters given in Figure~\ref{fig:silhouette graphs}.

The graphs in Figure~\ref{fig:silhouette graphs} indicate that one must
  consider the silhouette score and the \emph{marginal} silhouette scores simultaneously in
  order to determine the right number of clusters. For the synthetic datasets, the cophenetic metric produced the best silhouette scores. We also see that the silhouette score did indeed determine the right number of clusters for both of the synthetic datasets regardless of the choice of the metric. For the Iris datasets, all metrics appear to produce comparable silhouette scores. In this case, the cophenetic metric  suggests that the right number is 4 while other metrics agree on 2 clusters. We touch on
  this issue in the Conclusion. On the other hand, the consensus on the right number of
  clusters for the Cancer dataset is 2 for all metrics while cophenetic metric suggests 5.
  Finally, there appears to be no meaningful clusters for the Turkish cities dataset.

\subsubsection{Linkages}

There is one more parameter that effects the forming of clusters: the linkages we use to
calculate the distances between newly merged clusters as discussed in
Section~\ref{sect:linkages}. For this set of experiments, we used hierarchical clustering on
each dataset using different linkage methods.

Since hierarchical clustering is an unsupervised method, to use measures such as the
F1-score and accuracy, we need to find a suitable permutation of the confusion matrix for
each dataset. To deal with this problem, we use the Hungarian assignment method as defined
in~\cite{hungarian2005}. We then evaluate each model using F1-score (F1), accuracy (Acc.),
homogeneity (Hom.), completeness (Comp.)~\cite{homegcomplete}, mutual information
(M.Info)~\cite{mutual} and Rand index (Rand)~\cite{randscore}.  We display our results in
Table~\ref{table:alldataset}.  In each table, we report the best performing linkage method
(single (S), complete (C), average (A), or ward (W)) with each evaluation measure. We excluded the linearly separable synthetic dataset because all of the
  metrics did produce the most optimal result with the appropriate linkage.

\begin{table} 
	\caption{A comparison of metrics on  datasets.}\label{table:alldataset}
	\centering
	\begin{subtable}[h]{\textwidth}
		\centering	
		\caption{Iris dataset.}
		\footnotesize
		\begin{tabular}{ p{4\width} p{0.9\width} p{0.9\width} p{0.9\width} p{0.9\width} p{0.9\width} p{0.9\width} p{0.9\width} p{0.9\width} p{0.9\width} p{0.9\width} p{0.9\width} p{0.9\width}  }
			{\bf Metric} &\multicolumn{2}{ c }{{\bf F1}}
			&\multicolumn{2}{ c }{{\bf Acc.}}
			&\multicolumn{2}{ c }{{\bf Hom.}}
			&\multicolumn{2}{ c }{{\bf Comp.}}
			&\multicolumn{2}{ c }{{\bf M.Info}}
			&\multicolumn{2}{c}{{\bf Rand}}\\\hline
			bray-curtis     & 0.82 & W & 0.88 & W & 0.69 & W & 0.95 & S   & 0.75 & S & 0.82 & W  \\\hline
			cosine          & 0.68 & W & 0.79 & W & 0.58 & W & 0.95 & S   & 0.64 & S & 0.77 & W  \\\hline
			manhattan       & 0.88 & W & 0.92 & W & 0.74 & W & 0.92 & S   & 0.82 & A & 0.87 & W  \\\hline
			euclidean       & 0.88 & A & 0.92 & A & 0.77 & A & 0.95 & S   & 0.84 & A & 0.87 & A  \\\hline
			minkowski       & 0.85 & A & 0.90 & A & 0.77 & A & 0.95 & S   & 0.80 & A & 0.85 & A  \\\hline
			homological     & 0.76 & W & 0.84 & W & 0.58 & S & 1.00 & S   & 0.64 & S & 0.78 & S  \\\hline
		\end{tabular}
	\end{subtable}
	\newline\vspace*{0.2cm}
	
	\begin{subtable}[h]{\textwidth}
		\centering		
		\caption{Cancer Coimbra dataset.}
		\footnotesize
		\begin{tabular}{ p{4\width} p{0.9\width} p{0.9\width} p{0.9\width} p{0.9\width} p{0.9\width} p{0.9\width} p{0.9\width} p{0.9\width} p{0.9\width} p{0.9\width} p{0.9\width} p{0.9\width}  }
			{\bf Metric} &\multicolumn{2}{ c }{{\bf F1}}
			&\multicolumn{2}{ c }{{\bf Acc.}}
			&\multicolumn{2}{ c }{{\bf Hom.}}
			&\multicolumn{2}{ c }{{\bf Comp.}}
			&\multicolumn{2}{ c }{{\bf M.Info}}
			&\multicolumn{2}{c}{{\bf Rand}}\\\hline
			bray-curtis  & 0.56 & S  & 0.56 & S & 0.02 & A  & 0.14  & S   & 0.02 & A  & 0.50 & S  \\\hline
			cosine            & 0.55 & C  & 0.55 & C & 0.01 & S  & 0.12  & S   & 0.01 & S  & 0.50 & C  \\\hline
			manhattan      & 0.53 & S	& 0.53 & S & 0.02 & A  & 0.13  & A   & 0.02 & A  & 0.50 & S  \\\hline
			euclidean      & 0.54 & W  & 0.54 & W & 0.02 & A  & 0.13  & A   & 0.02 & A  & 0.50 & W  \\\hline
			minkowski      & 0.53 & S  & 0.53 & S & 0.02 & A  & 0.13  & A   & 0.02 & A  & 0.50 & S  \\\hline
			homological  & 0.61 & W  & 0.61 & W & 0.03 & W  & 1.00  & S   & 0.02 & W  & 0.52 & W  \\\hline
		\end{tabular}
	\end{subtable}
	\newline\vspace*{0.2cm}

	\begin{subtable}[h]{\textwidth}
		\caption{Synthetic dataset with mixing.}
		\centering
		\footnotesize
		\begin{tabular}{ p{4\width} p{0.9\width} p{0.9\width} p{0.9\width} p{0.9\width} p{0.9\width} p{0.9\width} p{0.9\width} p{0.9\width} p{0.9\width} p{0.9\width} p{0.9\width} p{0.9\width}  }
			{\bf Metric} &\multicolumn{2}{ c }{{\bf F1}}
			&\multicolumn{2}{ c }{{\bf Acc.}}
			&\multicolumn{2}{ c }{{\bf Hom.}}
			&\multicolumn{2}{ c }{{\bf Comp.}}
			&\multicolumn{2}{ c }{{\bf M.Info}}
			&\multicolumn{2}{c}{{\bf Rand}}\\\hline
			bray-curtis    & 1.00  & A  & 1.00  & A  & 1.00  & A & 1.00 & A & 1.38 & A  & 1.00 & A \\\hline
			cosine         & 0.83  & A  & 0.91  & A  & 0.72  & C & 0.77 & S & 1.38 & C  & 0.87 & A \\\hline
			manhattan      & 1.00  & S  & 1.00  & S  & 1.00  & S & 1.00 & S & 0.99 & S  & 1.00 & S \\\hline
			euclidean      & 1.00  & A  & 1.00  & A  & 1.00  & A & 1.00 & A & 1.38 & A  & 1.00 & A \\\hline
			minkowski      & 1.00  & C  & 1.00  & C  & 1.00  & C & 1.00 & C & 1.38 & C  & 1.00 & C \\\hline
			homological    & 0.98  & A  & 0.99  & A  & 0.95  & A & 1.00 & S & 1.31 & A  & 0.98 & A \\\hline
		\end{tabular}
	\end{subtable}
	
\end{table}

The cophenetic distance did produce best results for the Cancer dataset
  across the board. For the synthetic dataset with mixing, the results for the cophenetic
  distance appear to be on par with the other metrics even though it produced the weakest
  results after the cosine distance for the Iris dataset. Our results indicate that
cophenetic distance does produce competitive results on measures such as the F1-score,
accuracy, homogeneity and the Rand index while it shines on measures such as completeness
and mutual information score consistently on all datasets.  The results indicate that
cophenetic distance does tend to produce complete clusters that show high average
inter-class dissimilarities.

\section{Conclusions and future work}\label{sect:conc_future}

\subsection{Conclusions}

We defined a non-archimedean metric, called the cophenetic metric, on persistent homology
classes.  We then used this metric to sketch rooted tree presentations for zeroth persistent
homology classes instead of sketching rooted trees on points in the data set. We note that having a non-archimedean metric persistent homology classes in all degrees allows one to visualize higher homology classes as dendrograms as well.

Since the zeroth homology classes naturally correspond to connected components of the
subspace from which our data set is sampled, one can now compare the
  results of hierarchical clustering schemes with different metrics on data points with the
results we obtain from the cophenetic distance on the zeroth homology.  To test the
soundness of our proposal, we did numerical experiments on the geographical coordinates of a
small sample cities of Turkey to compare the dendrograms coming the cophenetic metric on the
zeroth homology and the dendrograms of hierarchical clustering algorithms by varying metrics
in Section~\ref{subsect:TurkishCities}.

The results of our numerical experiments we outlined in Section~\ref{subsect:TurkishCities}
indicate that there is a statistically verifiable strong correlation between the dendrograms
coming from the cophenetic distance matrix $C_0(D)$ and the dendrograms coming from other
metrics.  The statistical evidence we collected supports our hypothesis that hierarchical
clustering and zeroth persistent homology together with the cophenetic metric yield
statistically verifiable commensurate topological information about the connected components
of the datasets we used in our analyses.

We also note that while hierarchical clustering algorithms exclusively rely on a metric
structure on the data cloud alone, persistent homology relies on the simplicial technology
to derive its results, and therefore, should be impervious to the
  underlying metric.  On the other hand, the Mantel test results in
  Section~\ref{mantelsection} indicate that cophenetic distance and Euclidean distances are
  most similar. This may come from the fact that Vietoris-Rips complex we used to calculate
  our homological invariants uses the Euclidean metric for its filtration structure.

On the practical side, one has to also check that the cophenetic distance in combination
with hierarchical clustering methods does create high quality clusters.  For this purpose,
we compared clusters coming from hierarchical clustering algorithms on a handful of datasets
with varying number of clusters using different metrics including our own cophenetic metric
on the zeroth persistent homology in Section~\ref{subsect:OtherDatasets}. Results we obtain indicate that cophenetic distance shines in creating
  high quality complete clusters where distinct clusters have low mutual information.

 We notice that even though the original dataset has 3 preset clusters,
  Figure~\ref{fig:silhouette graphs} suggests that the optimal number of clusters for the
  Iris dataset is 2 for all metrics while the cophenetic distance suggests that it is
  4. However, one must observe that the classes \emph{iris versicolor} and \emph{iris
    virginica} are
  intertwined~\cite{zbMATH01810276,iris4cluster-fig-9-10,iris4cluster-table-4}. Our
  computations appear to detect this phenomenon.  The cophenetic distance result suggests
  splitting extra subclusters along their intersection while results from other metrics
  suggest merging these clusters.

A similar phenomenon appears in the Cancer Coimbra dataset. In
\cite{cancer}, the authors use logistic regression, random forests, and support vector
machines to label data points as \emph{control} or \emph{patient} with specificity and
sensitivity in the high 80\%'s. Our results, on the other hand, were in the 60\%'s for both
specificity and sensitivity. However, we focus on using hierarchical clustering algorithms
with commonly-used metrics and our cophenetic metric to split the dataset into meaningful
clusters instead of labeling data points as \emph{control} or
\emph{patient}. The results we obtained in
  Section~\ref{subsect:OtherDatasets} indicate that the dataset contains homogeneous
meaningful subsets other than \emph{control} and \emph{patient}.

\subsection{Future work}

One can extend the results of this article in different directions. The first obvious avenue for extension is replacing the zeroth homology with higher persistent homology. As we noted above, dendrograms as cobordisms of 0-spheres are adequate in representing the
  relationships between zeroth persistent homology classes. For the higher homology classes
we would need to deal with higher cobordisms of $n$-spheres~\cite{Strong:Cobordism} if we
are to develop a similar theory for the $n$-th persistent homology.  For the first
persistent homology, the cobordisms are given by genus-$g$ Riemann surfaces with punctures.
Fortunately, there is a complete classification of such surfaces in
full~\cite{Donaldson:RiemannSurfaces}.  Unfortunately, for higher dimensional homology, the
cobordisms require higher dimensional manifolds with finitely many punctures for which there
is no classification exists.

The second avenue of extension we would consider is extending our result to data sets that
cannot be easily embedded in an affine space.  This is often the case when one deals with
categorical data that require different techniques than numerical
data~\cite{Agresti:IntroductionToCategoricalDataAnalysis}.  We have shown that provided one
can define a simplicial complex out of data sets whose features are purely or partially
categorical, the cophenetic homological distance would yield usable information about the
data set on par with hierarchical clustering.

\subsection*{Acknowledgements}
	The first author was supported by Research Fund Project Number TDK-2020-42698 of the
	Istanbul Technical University.

\bibliographystyle{plain}   
\bibliography{article}

\end{document}